\documentclass[letterpaper, 10 pt, conference]{ieeeconf}

\IEEEoverridecommandlockouts    
\usepackage{textcomp}

\usepackage{sansmath}
\usepackage{graphicx}
\usepackage{amsmath}
\usepackage{xspace}
\usepackage{multirow}
\usepackage{diagbox}
\usepackage{amsfonts}
\usepackage{mathtools}
\usepackage{hyperref}
\usepackage[english]{babel}
\usepackage[version=4]{mhchem}
\usepackage{siunitx}
\usepackage{blindtext}
\usepackage{longtable}
\usepackage{tabularx}
\usepackage[utf8]{inputenc}
\usepackage{subcaption}
\usepackage{float}
\usepackage{titlesec}
\usepackage{tikz}
\usepackage{pgfplots}
\usepgfplotslibrary{groupplots}
\pgfplotsset{compat=1.17}
\usepackage{booktabs}
\usepackage{pgfplotstable}
\usepackage{amsthm}
\usepackage{bm}
\usepackage{cancel}
\usepackage{cite}
\usepackage{mathrsfs}
\usepackage{etoolbox}
\usepackage{amssymb}
\usepackage[ruled,vlined]{algorithm2e}
\newcommand{\argmin}{\mathop{\rm argmin}}
\newcommand{\argmax}{\mathop{\rm argmax}}

\newcommand{\eg}{{\it e.g.}}
\newcommand{\ie}{{\it i.e.}}

\DeclareMathOperator{\blkdiag}{blkdiag}

\newcommand{\mcl}[1]{\mathcal{#1}}
\newcommand{\mbb}[1]{\mathbb{#1}}
\newcommand{\mbf}[1]{\mathbf{#1}}
\newcommand{\mrm}[1]{\mathrm{#1}}

\DeclareMathOperator*{\minimize}{minimize}

\newcommand{\subjectto}{\mathop{\rm subject~to}}

\theoremstyle{definition}
\newtheorem{remark}{Remark}
\newtheorem{proposition}{Proposition}

\newtheorem{assumption}{Assumption}

\newcommand{\OURS}{ADMM-EKI}

\title{\LARGE \bf
Ensemble Kalman Inversion for Constrained Nonlinear MPC: \\ An ADMM-Splitting Approach
}

\author{Ahmed Khalil$^*$, Mohamed Safwat$^*$, Efstathios Bakolas
\thanks{* Denotes equal contribution.}
\thanks{Ahmed Khalil and Efstathios Bakolas are with the Department of Aerospace Engineering and Engineering Mechanics, University of Texas at Austin, Austin, TX 78712 USA. Mohamed Safwat is with the Department of Mechanical Engineering, University of Washington, Seattle, WA 98195 USA. (Emails: {\tt\small akhalil@utexas.edu; mohsaf@uw.edu; bakolas@austin.utexas.edu}).}
}

\begin{document}

\bstctlcite{IEEEReferenceControl}

\maketitle
\thispagestyle{empty}
\pagestyle{empty}

\begin{abstract}
This work proposes a novel Alternating Direction Method of Multipliers (ADMM)-based Ensemble Kalman Inversion (EKI) algorithm for solving constrained nonlinear model predictive control (NMPC) problems. First, stage-wise nonlinear inequality constraints in the NMPC problem are embedded via an augmented Lagrangian with nonnegative slack variables. We then show that the resulting unconstrained augmented-Lagrangian primal subproblem admits a Bayesian interpretation: under independent Gaussian virtual observations, its minimizers coincide with MAP estimators, enabling solution via EKI. However, since the nonnegativity constraint on the slacks is a hard constraint not naturally encoded by a Gaussian model, our proposed algorithm yields a two-block ADMM scheme that alternates between (i) an inexact primal step that minimizes the augmented-Lagrangian objective (implemented via EKI rollouts), (ii) a nonnegativity projection for the slacks, and (iii) a dual ascent step. To balance exploration and convergence, an annealing schedule tempers sampling covariances while a penalty schedule increases constraint enforcement over outer iterations, encouraging global search early and precise constraint satisfaction later. We evaluate the proposed controller on a 6-DOF UR5e manipulation benchmark in \texttt{MuJoCo}, comparing it against DIAL-MPC (an iterative MPPI variant) as the arm traverses a cluttered tabletop environment.
\end{abstract}

\section{Introduction}

Optimal control plays a central role in robotics, aerospace, and power systems, where it supports agile manipulation \cite{kuindersma2016atlas,neunert2018wbnmpc}, locomotion \cite{grandia2019feedback}, precision control \cite{safwat2026data}, fuel-efficient guidance and trajectory tracking \cite{blackmore2010minlanding}, and economic dispatch and grid stabilization \cite{low2014opf1,khalil2025DLCS}. In many of these applications, the resulting problems involve nonlinear dynamics and constraints, making real-time solutions challenging. This is the setting of nonlinear model predictive control (NMPC), which repeatedly solves a finite-horizon optimal control problem and applies only the first control input before replanning \cite{rawlings2020mpc,diehl2009nmpc,grune2017nmpc,Zinage2025TransformerMPC}.

A major class of NMPC methods is based on trajectory optimization via local second-order approximations, including Differential Dynamic Programming (DDP) \cite{jacobson1970ddp}, iLQR \cite{li2004ilqr}, and iLQG \cite{todorov2005ilqg}. These methods iteratively linearize the dynamics and quadratize the cost around a nominal trajectory, yielding efficient Riccati-style updates. They are highly effective when accurate derivatives are available, and the iterates remain within a favorable local region, but performance can degrade for strongly nonlinear or nonsmooth problems \cite{diehl2009nmpc}.

\begin{figure}[t]
\centering
\includegraphics[width=0.7\linewidth]{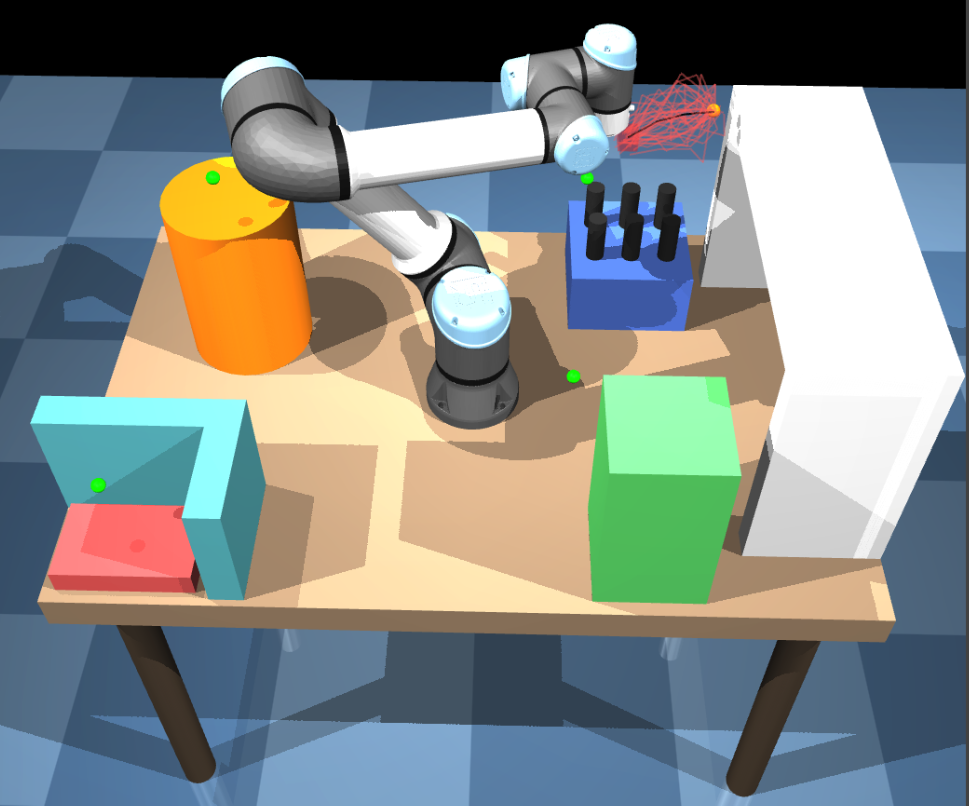}
\caption{UR5e robot running \OURS\ in a cluttered \texttt{MuJoCo} tabletop workcell while sequentially visiting predefined Cartesian waypoints (green spheres). Red lines show sampled rollout trajectories, and the black line denotes the corresponding mean trajectory.}
\label{fig_mujoco}
\end{figure}

A second major class is derivative-free sampling-based control, most notably Model Predictive Path Integral (MPPI) control \cite{kappen2005prl,theodorou2010jmlr,williams2017jgcd,williams2017icra}. MPPI updates control sequences by sampling trajectories, rolling out dynamics, and reweighting samples according to each trajectory's cost. Although a number of variants have been proposed to improve robustness and sample efficiency \cite{williams2018tubemppi,yin2022ccmppi,mohamed2022logmppi,yin2023shieldmppi,mohamed2023gpmppi,yan2024omppi}, principled handling of nonlinear constraints and nonconvex objectives remains challenging.

In this work, we propose a derivative-free constrained MPC method that combines Ensemble Kalman Inversion (EKI) with the Alternating Direction Method of Multipliers (ADMM). EKI is an ensemble-based derivative-free method for nonlinear inverse problems \cite{iglesias2013eki,schillings2017eki} whose use in optimal control is only beginning to emerge \cite{Askari2024EKS,joshi2024denkf}. Related iterative refinement ideas have recently appeared in sampling-based MPC, including DIAL-MPC \cite{xue2025dialmpc}, CoVO-MPC \cite{yi2024CoVOMPC}, and D4ORM \cite{zhang2025d4orm}. Our method uses EKI to approximately solve the ADMM primal subproblem, while ADMM handles nonlinear inequality constraints through a slack-variable split.

ADMM is a flexible and practically effective framework for constrained optimization \cite{boyd2011admm}. Compared with pure penalty methods \cite{carrillo2023consensus}, which can become ill-conditioned as penalty parameters grow \cite{nocedal2006numerical}, ADMM combines moderate penalties with multiplier updates to enforce constraints. In optimal control, related augmented-Lagrangian and ADMM ideas have shown strong performance in real-time settings, including TinyMPC \cite{nguyen2023tinympc} and ALTRO \cite{howell2019altro}.

The proposed method, \OURS, combines derivative-free rollouts, ensemble-based primal updates, and ADMM-based constraint handling. Its main features are:
\begin{itemize}
\item \textbf{Constrained nonlinear MPC:} stage-wise nonlinear inequality constraints are handled through an augmented-Lagrangian ADMM splitting.
\item \textbf{Derivative-free primal updates:} the primal subproblem is solved approximately using EKI and rollout-based sample covariances.
\item \textbf{Simple update structure:} each iteration consists of an EKI primal step, a slack projection, and a dual ascent step.
\item \textbf{Exploration and refinement:} annealing and penalty schedules encourage broader search early and tighter constraint enforcement later.
\end{itemize}

The remainder of this work is structured as follows. Section~\ref{section_problem_formulation} formulates the constrained NMPC problem. Section~\ref{section_ADMM_EKI} presents the proposed algorithm, \OURS. Section~\ref{section_numerical_simulations} evaluates \OURS\ in a 6-DOF UR5e manipulation benchmark in \texttt{MuJoCo} with obstacles, comparing against DIAL-MPC. Section~\ref{section_conclusion} concludes the paper and outlines future directions.

\section{Notation}
Let $\mbb{N}$ denote the set of nonnegative integers, $\mbb{R}^n$ the $n$-dimensional Euclidean space, and $\mbb{R}_+^n$ the nonnegative orthant. For a matrix $M \succeq 0$, define the weighted norm $\|v\|_M^2 \coloneqq v^\top M v$; $M \succ 0$ denotes that $M$ is positive definite. Inequalities between vectors are interpreted component-wise. For vectors $v_0,\ldots,v_H$ of compatible dimension, we use the stacked-vector convention $ [v_0^\top,\ldots,v_H^\top]^\top$. We write $I_p$ for the $p\times p$ identity matrix, $\otimes$ for the Kronecker product, and $\blkdiag(\cdot)$ for a block-diagonal matrix. The notation $[a]^+$ denotes the component-wise projection onto $\mbb{R}_+$, \ie, $[a]^+_j=\max\{a_j,0\}$. We use $\mcl{N}(\mu,\Sigma)$ to denote a Gaussian distribution with mean $\mu$ and covariance $\Sigma$. Throughout, $t$ indexes time, $\ell$ indexes ADMM outer iterations, $k$ indexes inner EKI iterations, and $i$ indexes ensemble particles.

\section{Problem Formulation}
\label{section_problem_formulation}

Let $H \in \mbb{N}$ be a finite prediction horizon and define the index set $ \mcl{H} \coloneqq \{0,\ldots,H-1\}$. Consider the (possibly nonlinear) discrete-time system
\begin{equation}
\label{eq_nonlinear_dynamics}
x_{t+1}=f(x_t,u_t), \qquad t \in \mcl{H},
\end{equation}
where $f:\mbb{R}^n \times \mbb{R}^m \to \mbb{R}^n$ is the one-step state-transition map, $x_t \in \mbb{R}^n$ is the state, $u_t \in \mbb{R}^m$ is the input, and the initial condition is known and given by $x_0=\bar{x} \in \mbb{R}^n$.

\begin{assumption}
\label{assumption_deterministic_markov}
For every initial condition $\bar{x}\in\mbb{R}^n$ and every input sequence $U \coloneqq [u_0^\top,\ldots,u_{H-1}^\top]^\top \in \mbb{R}^{Hm}$, the recursion \eqref{eq_nonlinear_dynamics} admits a unique solution /state trajectory $ X \coloneqq [x_0^\top,\ldots,x_H^\top]^\top \in \mbb{R}^{(H+1)n}$. 
\end{assumption}

Under Assumption~\ref{assumption_deterministic_markov}, the state trajectory is a deterministic function of the initial condition and the input sequence. We therefore define the trajectory rollout map $ \mcl{F}:\mbb{R}^n \times \mbb{R}^{Hm} \to \mbb{R}^{(H+1)n}$, by $ \mcl{F}(\bar{x},U) \coloneqq X$, where $X$ is the unique stacked trajectory generated by \eqref{eq_nonlinear_dynamics} from $(\bar{x},U)$.

Let $g:\mbb{R}^n \times \mbb{R}^m \to \mbb{R}^q$ be a stage-wise constraint map, interpreted component-wise, \ie, $g(x_t,u_t)\le 0$ means that each component of $g(x_t,u_t)$ is nonpositive. Additionally, define the stacked constraint map
\begin{align*}
\mcl{G}(X,U) \coloneqq
\begin{bmatrix}
g(x_0,u_0)\\
\vdots\\
g(x_{H-1},u_{H-1})
\end{bmatrix}
\in \mbb{R}^{Hq}.
\end{align*}

We consider the finite-horizon optimal control problem
\begin{equation}
\begin{aligned}
\label{eq_general_oct_problem}
&\minimize_{X,U} && \mcl{J}(X,U) \\
&\subjectto && X=\mcl{F}(\bar{x},U), \\
&&& \mcl{G}(X,U)\le 0,
\end{aligned}
\end{equation}
where $\mcl{J} : \mbb{R}^{(H+1)n} \times \mbb{R}^{Hm} \to \mbb{R}$ is the cost function.

\begin{assumption}
\label{assumption_cost_function}
The cost function $\mcl{J}(X,U)$ is quadratic.
\end{assumption}

In particular, for the remainder of the paper, we use the quadratic cost function
\begin{align*}
\mcl{J}(X, U) \coloneqq \tfrac{1}{2} \!\! \sum_{t=0}^{H-1} \left( \|x_t-z_t\|_R^2 {+} \|u_t\|_Q^2 \right) {+} \tfrac{1}{2} \|x_H-z_H\|_{R_H}^2,
\end{align*}
where $R \succ 0$, $R_H \succ 0$, and $Q \succ 0$ are weighting matrices, and $Z \coloneqq [z_0^\top,\ldots,z_H^\top]^\top$ is the reference trajectory.

\begin{remark}
Assumption~\ref{assumption_cost_function} does not render Problem~\eqref{eq_general_oct_problem} convex, since the dynamics constraint $X=\mcl{F}(\bar{x},U)$ and the inequality constraint $\mcl{G}(X,U)\le 0$ may still be nonlinear.
\end{remark}

\section{Ensemble Kalman Inversion using ADMM}
\label{section_ADMM_EKI}

In this section, we present the proposed \OURS\ algorithm. Subsection~\ref{subsec_augmented_lagrangian} introduces an augmented-Lagrangian treatment of the stage-wise inequality constraints via nonnegative slack variables. Subsection~\ref{subsec_admm_outer_loop} then casts the method as a two-block ADMM with primal, slack, and scaled-dual updates. Subsection~\ref{subsec_map_nmpc} provides a Bayesian interpretation of the primal subproblem, which allows us to solve the NMPC problem using estimation methods. Finally, Subsection~\ref{subsec_primal_update_eki} shows how the primal subproblem can be solved using EKI.

When needed, we write $U^{\ell,(k)}$ for clarity; in Algorithm~\ref{alg_admm_eki} we abbreviate $U^{(k)} \equiv U^{\ell,(k)}$ within a fixed $\ell$.

\subsection{Augmented Lagrangian formulation}
\label{subsec_augmented_lagrangian}

By Assumption~\ref{assumption_deterministic_markov}, the state trajectory is uniquely determined by $(\bar{x},U)$. We therefore introduce the reduced cost and reduced constraint maps:
\begin{align*}
\mcl{J}_{\bar{x}}(U) \coloneqq \mcl{J}(\mcl{F}(\bar{x},U),U), \qquad \mcl{G}_{\bar{x}}(U) \coloneqq \mcl{G}(\mcl{F}(\bar{x},U),U).
\end{align*}
The inequality constraint $\mcl{G}_{\bar{x}}(U)\le 0$ is equivalently written by introducing a nonnegative slack variable $S\in\mbb{R}_+^{Hq}$ such that $ \mcl{G}_{\bar{x}}(U)+S=0$. With multiplier $\Lambda\in\mbb{R}^{Hq}$ and penalty parameter $\rho>0$, the augmented Lagrangian is
\begin{equation*}
\mcl{L}_\rho(U,S,\Lambda) = \mcl{J}_{\bar{x}}(U) +\Lambda^\top\big(\mcl{G}_{\bar{x}}(U)+S\big) +\tfrac{\rho}{2}\|\mcl{G}_{\bar{x}}(U)+S\|_2^2.
\end{equation*}
Here $\Lambda$ is unrestricted, since it is the multiplier for the equality constraint $\mcl{G}_{\bar{x}}(U)+S=0$; nonnegativity is enforced through the slack constraint $S\ge 0$. Introducing the scaled dual variable $Y\coloneqq \Lambda/\rho$ and completing the square gives
\begin{equation*}
\mcl{L}_\rho(U,S,Y) = \mcl{J}_{\bar{x}}(U) +\tfrac{\rho}{2}\big\|\mcl{G}_{\bar{x}}(U)+S+Y\big\|_2^2 -\tfrac{\rho}{2}\|Y\|_2^2.
\end{equation*}
The constant term $-\tfrac{\rho}{2}\|Y\|_2^2$ is omitted below since it does not affect minimizers in $(U,S)$. Note that the dynamics constraint is enforced implicitly through the reduced maps $\mcl{J}_{\bar{x}}$ and $\mcl{G}_{\bar{x}}$ and is therefore not dualized.

\subsection{ADMM Outer Loop}
\label{subsec_admm_outer_loop}

The two-block ADMM iteration, with $U$ in one block and $S$ in the other, consists of the following steps.

The primal update solves
\begin{equation}
\label{eq_primal_step}
U^{\ell+1}=\argmin_U \Phi^\ell(U),
\end{equation}
where
\begin{equation} \label{eq_PhiU}
\Phi^\ell(U) \coloneqq \mcl{J}_{\bar{x}}(U) +\tfrac{\rho^\ell}{2}\big\|\mcl{G}_{\bar{x}}(U)+S^\ell+Y^\ell\big\|_2^2.
\end{equation}

The slack update solves
\begin{align*}
S^{\ell+1} &= \argmin_{S\ge 0} \tfrac{\rho^\ell}{2} \big\|\mcl{G}_{\bar{x}}(U^{\ell+1})+S+Y^\ell\big\|_2^2,
\end{align*}
which has the closed-form solution
\begin{equation} \label{eq_slack_update} S^{\ell+1} = \big[-\mcl{G}_{\bar{x}}(U^{\ell+1})-Y^\ell\big]^+.
\end{equation}

The scaled-dual update is
\begin{equation}
\label{eq_dual_step}
Y^{\ell+1} = Y^\ell+\mcl{G}_{\bar{x}}(U^{\ell+1})+S^{\ell+1}.
\end{equation}

The penalty parameter is updated as $\rho^{\ell+1}=\tau\rho^\ell$ with $\tau\ge 1$. Since $Y^\ell=\Lambda^\ell/\rho^\ell$, changing $\rho$ requires a corresponding rescaling of $Y$ to preserve this relationship. Thus, after updating $\rho$, we set
\begin{equation} \label{eq_dual_rescale} Y^{\ell+1} \leftarrow \frac{\rho^\ell}{\rho^{\ell+1}}\,Y^{\ell+1} = \frac{1}{\tau}Y^{\ell+1}.
\end{equation}

\subsection{Bayesian estimation view of the primal subproblem}
\label{subsec_map_nmpc}

We now give a Bayesian interpretation of \eqref{eq_primal_step}. Since $\bar{x}$ is fixed, the nonlinear map appearing in the ADMM penalty term is precisely the reduced constraint map $\mcl{G}_{\bar{x}}(U)$. Define
\begin{align*}
y^\ell \coloneqq -S^\ell-Y^\ell \in \mbb{R}^{Hq}.
\end{align*}
Additionally, define $\Sigma_R \coloneqq \blkdiag(I_H\otimes R,\ R_H)$, $\Sigma_Q \coloneqq I_H \otimes Q$, and $ \Sigma_{\rho^\ell}\coloneqq \frac{1}{\rho^\ell} I_{Hq}$.

The quantities $0$, $Z$, and $y^\ell$ are interpreted as \emph{virtual observations}, meaning auxiliary observations introduced for the Bayesian reformulation rather than physical measurements. In particular, $y^\ell$ is an ADMM-generated pseudo-datum for the reduced constraint map $\mcl{G}_{\bar{x}}(U)$, chosen so that the resulting negative log-likelihood reproduces the quadratic penalty term in \eqref{eq_PhiU}.

\begin{assumption}[Independent Gaussian virtual observations]
\label{assumption_observations}
Assume three independent Gaussian virtual observations:
\begin{align}
0 &= U + \eta_U, & \eta_U &\sim \mcl{N}(0,\Sigma_Q), \nonumber\\
Z &= \mcl{F}(\bar{x},U) + \eta_X, & \eta_X &\sim \mcl{N}(0,\Sigma_R), \nonumber\\
y^\ell &= \mcl{G}_{\bar{x}}(U) + \eta_G, & \eta_G &\sim \mcl{N}(0,\Sigma_{\rho^\ell}). \label{eq_virtual_obs}
\end{align}
\end{assumption}

\begin{proposition}
\label{prop_map_equals_min}
Let $\Phi^\ell(U)$ be as defined in \eqref{eq_PhiU}, and assume that the set of minimizers of $\Phi^\ell$ is nonempty. Under Assumption~\ref{assumption_observations}, the set of MAP estimators for the posterior $p(U\mid 0,Z,y^\ell)$ coincides with the set of minimizers of $\Phi^\ell$:
\begin{equation*}
\argmax_U p(U\mid 0,Z,y^\ell) = \argmin_U \Phi^\ell(U).
\end{equation*}
\end{proposition}

\begin{proof}
By Assumption~\ref{assumption_observations}, the conditional densities are Gaussian and independent given $U$, so
\begin{equation*}
p(U\mid 0,Z,y^\ell) \propto p(0\mid U)\,p(Z\mid U)\,p(y^\ell\mid U).
\end{equation*}
Moreover,
\begin{align*}
p(0\mid U) &\propto \exp\!\left(-\tfrac{1}{2}\|U\|_{\Sigma_Q^{-1}}^2\right),\\ 
p(Z\mid U) &\propto \exp\!\left(-\tfrac{1}{2}\|\mcl{F}(\bar{x},U)-Z\|_{\Sigma_R^{-1}}^2\right),\\
p(y^\ell\mid U) &\propto \exp\!\left(-\tfrac{1}{2}\|\mcl{G}_{\bar{x}}(U)-y^\ell\|_{\Sigma_{\rho^\ell}^{-1}}^2\right).
\end{align*}
Taking negative logarithms and discarding constants yields
\begin{multline*}
-\log p(U\mid 0,Z,y^\ell) = \\
\tfrac{1}{2} \|U\|_{\Sigma_Q^{-1}}^2  + \tfrac{1}{2} \|\mcl{F}(\bar{x},U)-Z\|_{\Sigma_R^{-1}}^2 + \tfrac{1}{2} \|\mcl{G}_{\bar{x}}(U)-y^\ell \|_{\Sigma_{\rho^\ell}^{-1}}^2 \\
= \tfrac{1}{2} \|U\|_{\Sigma_Q^{-1}}^2 + \tfrac{1}{2} \|\mcl{F}(\bar{x},U)-Z\|_{\Sigma_R^{-1}}^2 + \tfrac{\rho^\ell}{2} \| \mcl{G}_{\bar{x}}(U)-y^\ell \|_2^2.
\end{multline*}
By the quadratic cost definition and $y^\ell=-S^\ell-Y^\ell$,
\begin{align*}
\mcl{J}_{\bar{x}}(U) = \tfrac{1}{2} \|U\|_{\Sigma_Q^{-1}}^2 + \tfrac{1}{2} \|\mcl{F}(\bar{x},U)-Z\|_{\Sigma_R^{-1}}^2.
\end{align*}
Therefore,
\begin{align*}
-\log p(U\mid 0,Z,y^\ell) &= \mcl{J}_{\bar{x}}(U)  + \tfrac{\rho^\ell}{2}\|\mcl{G}_{\bar{x}}(U) + S^\ell+Y^\ell\|_2^2 \\ 
&= \Phi^\ell(U),
\end{align*}
which proves the claim.
\end{proof}

\begin{remark}
Since $\mcl{F}$ and $\mcl{G}_{\bar{x}}$ are generally nonlinear, the primal subproblem \eqref{eq_primal_step} is generally nonconvex. We therefore solve \eqref{eq_primal_step} inexactly, \ie, EKI is used to compute an approximate minimizer rather than the exact solution of the primal subproblem.
\end{remark}

\subsection{Primal update via Ensemble Kalman Inversion}
\label{subsec_primal_update_eki}

At each ADMM outer iteration $\ell$, we solve the primal subproblem in Equation \eqref{eq_primal_step} using Ensemble Kalman Inversion (EKI). The key idea is to maintain an ensemble of candidate control sequences and iteratively move that ensemble toward lower values of the primal objective. For each particle, we evaluate three quantities: the control sequence itself, the tracking error induced by the rollout $\mcl{F}(\bar{x},U)$, and the current ADMM constraint mismatch $\mcl{G}_{\bar{x}}(U)+S^\ell+Y^\ell$. These are stacked into a residual vector. EKI then uses empirical covariances between the ensemble of controls and the ensemble of residuals to construct a Kalman-type gain, which updates each particle in a direction that reduces these residuals. After $M$ inner EKI iterations, the ensemble mean is taken as the approximate primal iterate $U^{\ell+1}$.

For particle $i$ at inner iteration $k$, define the residual vector
\begin{equation}
\label{eq_residual}
C_i^{(k)}(\bar{x},U_i^{(k)},S^\ell,Y^\ell) \coloneqq
\begin{bmatrix}
U_i^{(k)} \\
\mcl{F}(\bar{x},U_i^{(k)})-Z \\
\mcl{G}_{\bar{x}}(U_i^{(k)})+S^\ell+Y^\ell
\end{bmatrix}
\in \mbb{R}^d,
\end{equation}
where $d \coloneqq Hm + (H+1)n + Hq$. The first block corresponds to the control regularization term, the second to the trajectory tracking error, and the third to the augmented-Lagrangian penalty term.

Define the block-diagonal weighting matrix
\begin{equation}
\label{eq_q_hat}
\widehat{Q}_{\rho^\ell} \coloneqq \blkdiag(\Sigma_Q,\Sigma_R,\Sigma_{\rho^\ell}) \in \mbb{R}^{d\times d}.
\end{equation}

Initialize particles around a nominal mean $\bar{U}^{(0)}$ given by
\begin{equation*}
U_i^{(0)}=\bar{U}^{(0)}+\epsilon_i^{(0)}, \qquad \epsilon_i^{(0)} \sim \mcl{N}(0,\beta^{(0)}\Sigma_U),
\end{equation*}
where $\Sigma_U$ is a design covariance and $\beta^{(k)}$ is an annealing parameter that controls the amount of injected exploration.

Define the ensemble means
\begin{equation} \label{eq_ensemble_means}
\bar{U}^{(k)}=\frac{1}{N}\sum_{i=1}^N U_i^{(k)}, \qquad \bar C^{(k)}=\frac{1}{N}\sum_{i=1}^N C_i^{(k)}.
\end{equation}
Define the anomaly matrices
\begin{align*}
\Delta U^{(k)} &= \big[ U_1^{(k)}-\bar{U}^{(k)},\ldots,U_N^{(k)}-\bar{U}^{(k)} \big] \in \mbb{R}^{Hm\times N}, \\
\Delta C^{(k)} &= \big[ C_1^{(k)}-\bar C^{(k)},\ldots,C_N^{(k)}-\bar C^{(k)} \big] \in \mbb{R}^{d\times N}.
\end{align*}
The corresponding sample covariances are
\begin{align}
P_{UC}^{(k)} &= \tfrac{1}{N-1}\Delta U^{(k)}(\Delta C^{(k)})^\top \in \mbb{R}^{Hm\times d}, \label{eq_cross_cov} \\
P_{CC}^{(k)} &= \tfrac{1}{N-1}\Delta C^{(k)}(\Delta C^{(k)})^\top \in \mbb{R}^{d\times d}. \label{eq_CC_cov}
\end{align}

Using these empirical covariances, we compute the ensemble Kalman gain
\begin{equation}
\label{eq_EKI_gain_update}
K^{(k)} = P_{UC}^{(k)} \Big(P_{CC}^{(k)}+\widehat{Q}_{\rho^\ell}\Big)^{-1}.
\end{equation}
Each particle is then updated by
\begin{equation} \label{eq_EKI_update_law}
U_i^{(k+1)} = U_i^{(k)} - K^{(k)}C_i^{(k)} + \epsilon_i^{(k+1)},
\end{equation}
where the additive jitter term is given by $ \epsilon_i^{(k+1)}\sim\mcl{N}(0,\beta^{(k+1)}\Sigma_U)$. This perturbation promotes exploration in the early iterations and gradually vanishes as the annealing parameter decreases.

We use the exponential decay schedule
\begin{equation} \label{eq_annealing}
\beta^{(k+1)}=\beta_0 e^{-\gamma(k+1)},
\end{equation}
with $\beta_0>0$ and $\gamma>0$.

\begin{algorithm}[h]
\caption{\OURS\ Model Predictive Control}
\label{alg_admm_eki}
\SetAlgoLined
\LinesNumbered
\SetKwInOut{KwIn}{Input}

\KwIn{Parameters $\Sigma_U, N, L, M, \tau, \rho^0, \beta_0, \gamma$}
\KwIn{Warm start $(S^0,Y^0,\bar{U}^{(0)})$}

\While{task not complete}{
$\bar{x} \leftarrow \textit{GetStateEstimate}()$ \\
\For{$\ell=0,1,\ldots,L$}{
\textbf{[Weights]} $\widehat{Q}_{\rho^\ell}\gets \blkdiag(\Sigma_Q,\Sigma_R,\Sigma_{\rho^\ell})$ \\
\textbf{[Ensemble]} \For{$i=1,\ldots,N$ \textbf{in parallel}}{
$\epsilon_i^{(0)}\sim \mcl{N}(0,\beta^{(0)}\Sigma_U)$ \\
$U_i^{(0)}\gets \bar{U}^{(0)}+\epsilon_i^{(0)}$
}
\textbf{[Mean]} $\bar{U}^{(0)}\gets \frac{1}{N}\sum_{i=1}^N U_i^{(0)}$ \\
\For{$k=0,1,\ldots,M-1$}{
\For{$i=1,\ldots,N$ \textbf{in parallel}}{
\textbf{[Rollout]} $C_i^{(k)}(\bar{x},U_i^{(k)},S^\ell,Y^\ell)\gets$ \eqref{eq_residual}
}
\textbf{[Means]} $\bar{U}^{(k)},\bar C^{(k)}\gets$ \eqref{eq_ensemble_means} \\
\textbf{[Covariances]} $P_{UC}^{(k)},P_{CC}^{(k)}\gets$ \eqref{eq_cross_cov}, \eqref{eq_CC_cov} \\
\textbf{[Gain]} $K^{(k)}\gets$ \eqref{eq_EKI_gain_update} \\
\textbf{[Anneal]} $\beta^{(k+1)}\gets \beta_0 e^{-\gamma(k+1)}$ \\
\For{$i=1,\ldots,N$ \textbf{in parallel}}{
$\epsilon_i^{(k+1)}\sim \mcl{N}(0,\beta^{(k+1)}\Sigma_U)$ \\
\textbf{[EKI step]} $U_i^{(k+1)}\gets U_i^{(k)} - K^{(k)}C_i^{(k)} + \epsilon_i^{(k+1)}$
}
\textbf{[Mean]} $\bar{U}^{(k+1)}\gets \frac{1}{N}\sum_{i=1}^N U_i^{(k+1)}$ \\
}
\textbf{[Primal]} $U^{\ell+1}\gets \bar{U}^{(M)}$ \\
\textbf{[Slack]} $S^{\ell+1}\gets\big[-\mcl{G}_{\bar{x}}(U^{\ell+1})-Y^\ell\big]^+$ \\
\textbf{[Dual]} $Y^{\ell+1}\gets Y^\ell+\mcl{G}_{\bar{x}}(U^{\ell+1})+S^{\ell+1}$ \\
\textbf{[Penalty]} $\rho^{\ell+1}\gets \tau\rho^\ell$ \\
\textbf{[Rescale]} $Y^{\ell+1}\gets (\rho^\ell/\rho^{\ell+1})Y^{\ell+1}$
}
$\textit{ExecuteCommand}((U^{L+1})_0)$
}
\end{algorithm}

\begin{remark}[Interpretation of the EKI step]
The update \eqref{eq_EKI_update_law} moves each control particle according to the empirical correlation between variations in $U$ and variations in the stacked residual $C$. Consequently, particles are steered toward control sequences that simultaneously reduce control effort, improve trajectory tracking, and decrease the ADMM constraint mismatch. The final ensemble mean after $M$ inner iterations is used as the approximate solution of the primal subproblem.
\end{remark}

\begin{remark}[Efficient inversion via the Woodbury identity]
When $d\gg N$, the inverse in \eqref{eq_EKI_gain_update} can be computed efficiently using the matrix inversion lemma. Let $S\coloneqq \widehat{Q}_{\rho^\ell}$. Then, by \eqref{eq_CC_cov},
\begin{multline*}
\big(P_{CC}^{(k)}+S\big)^{-1} = S^{-1} - S^{-1}\Delta C^{(k)} \\
\cdot \Big( (N-1)I+(\Delta C^{(k)})^\top S^{-1}\Delta C^{(k)} \Big)^{-1} (\Delta C^{(k)})^\top S^{-1}.
\end{multline*}
This reduces a $d\times d$ inversion to an $N\times N$ inversion.
\end{remark}

Algorithm~\ref{alg_admm_eki} summarizes the resulting method. At a high level, each MPC call begins from the current state estimate $\bar{x}$ and a warm start from the previous timestep. For each ADMM outer iteration, an ensemble of candidate control sequences is initialized around the warm start and evolved through $M$ EKI updates. These updates use repeated rollouts of the nonlinear dynamics to evaluate residuals and empirical covariances, which in turn define the Kalman gain. After the inner loop terminates, the ensemble mean is taken as the new primal variable. The slack and scaled dual variables are then updated using the standard ADMM steps, the penalty parameter is increased if desired, and the procedure repeats. After the final outer iteration, only the first control input of the optimized sequence is applied, and the remaining quantities are used to warm-start the next receding-horizon solve.

\subsection{Receding-horizon implementation of \OURS}

We run Algorithm~\ref{alg_admm_eki} in a receding-horizon MPC loop, executing only the first control input of the optimized sequence. Each call is warm-started with $(S^0,Y^0,\bar{U}^{(0)})$ from the previous timestep, or initialized to zero at the first timestep. Given the current state estimate, we perform $L+1$ ADMM outer iterations, each containing $M$ EKI updates, then apply the first element of the final mean control sequence and repeat at the next sampling instant.

\begin{remark}
\label{rem_termination}
Rather than always performing $L+1$ outer iterations, one may terminate early using standard ADMM primal and dual residual tests and adapt $\rho$ accordingly; see \cite{boyd2011admm}. We use fixed schedules here for simplicity.
\end{remark}

\section{Numerical Simulations} \label{section_numerical_simulations}

We evaluate the proposed controller, \OURS, on a 6-DOF UR5e waypoint-reaching benchmark in a cluttered \texttt{MuJoCo} tabletop workcell, shown in Fig.~\ref{fig_mujoco}. Each episode consists of sequentially reaching a list of Cartesian waypoints $\{\mbf{g}_i\}_{i=1}^{N_{\mrm{wp}}}\subset\mbb{R}^3$ in end-effector (EEF) position space, with $N_{\mrm{wp}}=5$ waypoints per episode. The objective is to visit all waypoints in order while avoiding collisions with obstacles and preventing self-collisions. The Cartesian waypoints are deliberately placed in close proximity to obstacles, stressing the controller's ability to reach each target while maintaining safe clearance throughout the motion. We run $10$ episodes per controller, and to probe robustness, the waypoint set is perturbed at the start of each episode as $\tilde{\mbf{g}}_i=\mbf{g}_i+\boldsymbol{\eta}_i$, where $\boldsymbol{\eta}_i\sim\mcl{N}(\mbf{0},\sigma^2 I_3)$ and $\sigma=0.02~\mrm{m}$.

The robot optimizes sequences of joint increments $\Delta \mbf{q}_t \in \mbb{R}^6$ over a horizon $H$ (\ie, $u_t \coloneqq \Delta \mbf{q}_t$). Given the current joint configuration $\mbf{q}_0$, candidate sequences are rolled out using the deterministic kinematic recursion
\begin{align*}
\mbf{q}_{t+1} = \mbf{q}_t + \Delta \mbf{q}_t,\qquad t \in \mcl{H},
\end{align*}
with elementwise clipping $\|\Delta \mbf{q}_t\|_\infty \le \Delta q_{\max}$ and $\Delta q_{\max}=0.12~\mrm{rad}$. For each rollout, we compute EEF positions $\mbf{x}_{\mrm{eef}}(t)=\mrm{FK}(\mbf{q}_t)$, where $\mrm{FK}:\mbb{R}^6\rightarrow \mbb{R}^3$ denotes the (\texttt{MuJoCo}-evaluated) forward-kinematics map from joint configuration to EEF position in the world frame. We also compute the relevant constraint signals (\eg, minimum signed distance to obstacles/robot collision geometries and joint-limit margins) directly from \texttt{MuJoCo}. In these experiments, the reference trajectory $Z$ in the quadratic cost (Assumption~\ref{assumption_cost_function}) is chosen to be terminal-only: we penalize only the final EEF position error at the end of the horizon. This yields a smaller, structured EKI residual (the $\Sigma_R^{-1}$ block acts only at $t{=}H$), reducing gain-update overhead.

We compare against DIAL-MPC~\cite{xue2025dialmpc}, an iterative MPPI-style sampling MPC method. At each control step, DIAL-MPC maintains a mean increment sequence $\mbf{U}\in\mbb{R}^{H\times 6}$ and draws $K$ Gaussian-perturbed candidates. Each sample is rolled out kinematically and scored by a cost comprising (i) EEF waypoint tracking, (ii) control effort, and (iii) soft penalties for collision proximity/penetration and joint-limit violation. DIAL-MPC updates the mean sequence via MPPI-style exponential reweighting with temperature $\lambda$. We use $H=5$, $K=40$, $\lambda=3.0$, and $5$ internal annealing stages per MPC step.

While MPPI-style methods can incorporate ``hard'' constraints by rejecting infeasible samples or using indicator-style costs, in cluttered, high-dimensional settings this often results in few or no feasible rollouts, weight degeneracy, and unstable updates unless the sampling distribution is carefully tuned \cite{yin2022ccmppi}. In contrast to relying on soft penalty terms in the cost, \OURS\ enforces collision avoidance and joint limits as explicit stage-wise inequality constraints via an augmented-Lagrangian split, thereby driving feasibility through slack projection and dual updates. For both controllers, we use the same per-stage constraint signals $\mbf{g}(t) = [g_{\mrm{coll}}(t),\, g_{\mrm{jl}}(t)]^\top$, where $g_{\mrm{coll}}(t)=d_{\mrm{safe}}-d_{\min}(t)$ encodes a minimum-distance safety margin and $g_{\mrm{jl}}(t)$ encodes joint-limit violation, stacked over the horizon as $\mcl{G}(\bar{x},U)$. Here $d_{\min}(t)$ is the minimum signed distance returned by \texttt{MuJoCo} over all relevant contact pairs, including robot--obstacle and robot self-contact pairs; in our experiments, we set $d_{\mrm{safe}}=0$~m (\ie, enforcing non-penetration, allowing at most contact at zero distance). Joint limits are taken from the UR5e model bounds in \texttt{MuJoCo} (the per-joint $[q_{\min},q_{\max}]$ ranges), and we use a zero margin in this benchmark. We use $H=5$, $N=40$ particles, $\rho^0=10$, and $\tau=1.8$. In Algorithm~\ref{alg_admm_eki}, we set $L=1$ and $M=5$ for each MPC call.

\begin{figure}[t]
\centering
\includegraphics[width=\linewidth]{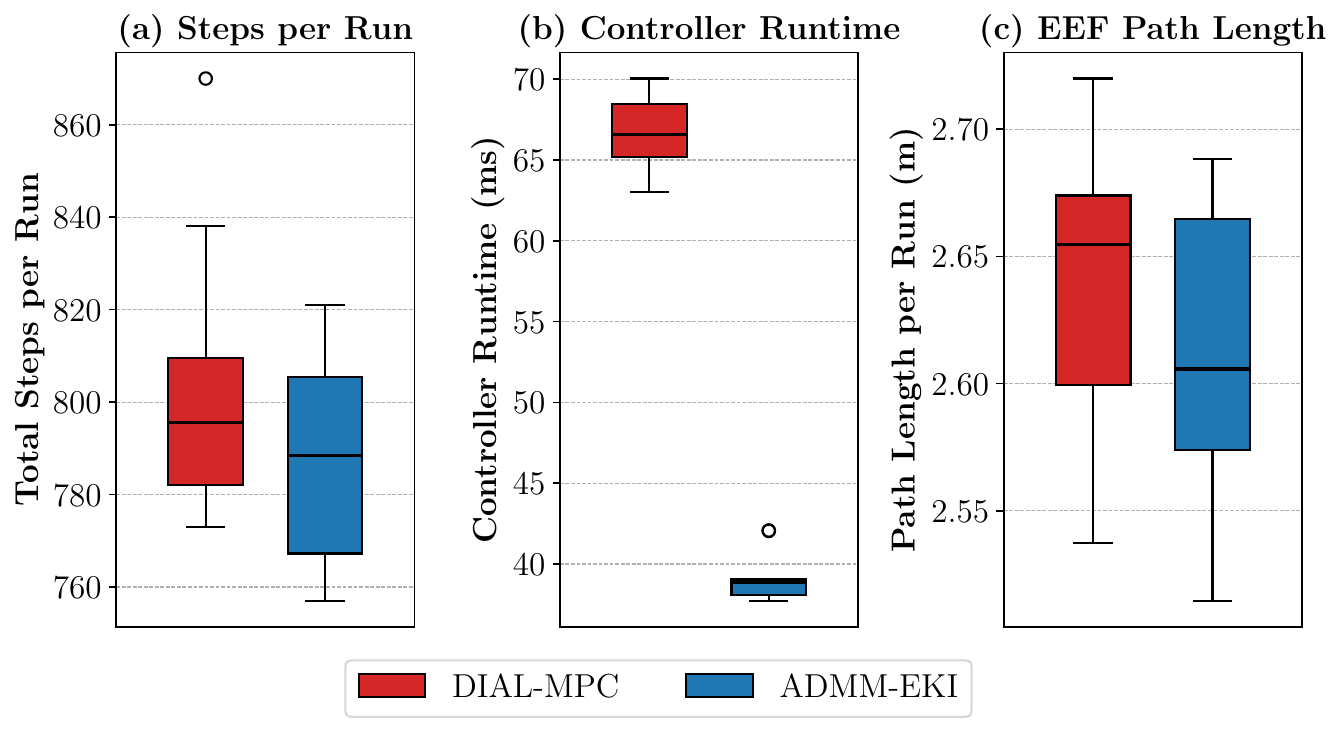}
\caption{Box plots of (a) total steps, (b) controller runtime (ms) per MPC update, and (c) end-effector path length across $10$ episodes, comparing DIAL-MPC and \OURS\ on the obstacle-cluttered UR5e waypoint task (lower is better).}
\label{fig_benchmarks}
\end{figure}

For each controller and episode, we record: (1) the total number of control steps (\ie, MPC/environment steps) until termination or the step budget is reached; (2) the controller runtime per MPC update in milliseconds; and (3) the end-effector path length over the episode, computed as $\sum_t \|\mbf{x}_{\mrm{eef}}(t{+}1)-\mbf{x}_{\mrm{eef}}(t)\|_2$ (lower is better for all metrics).

All benchmarking experiments were conducted on a laptop with an AMD Ryzen 9 4900HS CPU (3.0\,GHz) and integrated Radeon graphics. We summarize results across $10$ episodes with box plots (Fig.~\ref{fig_benchmarks}). Overall, \OURS\ outperforms DIAL-MPC on all metrics and is substantially faster: it requires fewer control steps (mean $787.2$ vs.\ $802.7$; max $821$ vs.\ $870$), reduces per-update runtime by $\approx41\%$ (mean $39.18$~ms vs.\ $66.75$~ms; $\approx1.7\times$ faster), and yields a shorter end-effector path length (mean $2.613$~m vs.\ $2.637$~m).

\section{Conclusion} \label{section_conclusion}

This work introduced \OURS, a derivative-free and parallelizable NMPC method that combines ensemble Kalman inversion with an ADMM-style augmented Lagrangian to handle stage-wise nonlinear inequality constraints, using annealing to balance exploration and refinement. On a 6-DOF UR5e waypoint-reaching benchmark in \texttt{MuJoCo}, \OURS\ achieves substantially lower per-step runtime than DIAL-MPC while improving steps-to-completion and end-effector path length. Future work includes adaptive ADMM penalty/termination strategies.

\bibliographystyle{IEEEtran}
\bibliography{root}

\end{document}